\documentclass[a4paper,12pt]{article}
\usepackage{amsmath}
\usepackage{amssymb}
\usepackage{amsfonts}
\begin{document}
\newtheorem{theoremenv}{Theorem}
\newtheorem{propositionenv}[theoremenv]{Proposition}
\newtheorem{lemmaenv}[theoremenv]{Lemma}
\newtheorem{corollaryenv}[theoremenv]{Corollary}
\newtheorem{claimenv}[theoremenv]{Claim}
\newtheorem{constructionenv}[theoremenv]{Construction}
\newtheorem{conjectureenv}[theoremenv]{Conjecture}
\newtheorem{questionenv}[theoremenv]{Question}
\newtheorem{definitionenv}[theoremenv]{Definition}
\renewcommand{\labelenumi}{(\roman{enumi})}
\newsavebox{\proofbox}
\savebox{\proofbox}{\begin{picture}(7,7)%
  \put(0,0){\framebox(7,7){}}\end{picture}}

\def\proof{\noindent\textbf{Proof. }}
\def\endproof{\hfill{\usebox{\proofbox}}}
\def\remark{\noindent\textbf{Remark. }}
\def\scs{\scriptstyle}

\newcommand{\Hom}{\mbox{Hom}}
\newcommand{\Cov}{\mbox{Cov}}
\newcommand{\Var}{\mbox{Var}}
\newcommand{\Supp}{\mbox{Supp}}
\newcommand{\Diam}{\mbox{Diam}}
\newcommand{\Span}{\mbox{Span}}
\newcommand{\Prob}{\mbox{Prob}}
\newcommand{\RanProb}{\mathbb{P}_{\mbox{ran}}}
\newcommand{\SetProb}{\mathbb{P}_{\mbox{set}}}
\newcommand{\RegProb}{\mathbb{P}_{\mbox{reg}}}

\begin{center}
\Large The Cameron--Erd\H{o}s Conjecture\\[15pt]
\large Ben Green\footnote{Supported by a Fellowship of Trinity College, Cambridge and a grant from the EPSRC, United Kingdom. Mathematics Subject Classification: 11B75.}
\end{center}
\begin{abstract} \noindent A subset $A$ of the integers is said to be \textit{sum-free} if there do not exist elements $x,y,z \in A$ with $x + y = z$. It is shown that the number of sum-free subsets of $\{1,\dots,N\}$ is $O(2^{N/2})$, confirming a well-known conjecture of Cameron and Erd\H{o}s.
\end{abstract}
\noindent\textbf{1. Introduction.} If $A$ is any subset of an abelian group then we say that $A$ is sum-free if $(A + A) \cap A = \emptyset$, that is if there do not exist $x,y,z \in A$ for which $x + y = z$. The study of such sets goes back at least 30 years, and over 10 years ago Cameron and Erd\H{o}s \cite{CE1,CE2} raised the question of enumerating the sum-free subsets of $[N] = \{1,\dots,N\}$. They noted that any set of odd integers is sum-free, as is any subset of $\{\lceil N/2 \rceil,\dots, N\}$, but that it is hard to think of many sum-free sets which are not essentially of this form. Thus they advanced the following conjecture.
\begin{conjectureenv}[Cameron--Erd\H{o}s] \label{conj1}The number of sum-free subsets of $[N]$ is $O(2^{N/2})$.
\end{conjectureenv}
There has been some progress on this conjecture. Writing $\mbox{SF}(N)$ for the collection of sum-free subsets of $[N]$, Alon \cite{Alon}, Calkin \cite{Calk} and Erd\H{o}s and Granville (unpublished) showed independently that
\begin{equation}\label{ACEG}
|\mbox{SF}(N)| \; = \; 2^{N/2 + o(N)}.
\end{equation}
Results in a rather different direction were obtained by Freiman \cite{Frei4} and by Deshouillers, Freiman, S\'os and Temkin \cite{DFST}. In \cite{Frei4}, for example, it was shown that the number of sum-free subsets of $[N]$ with cardinality at least $5N/12 + 2$ is at most $O(2^{N/2})$. Let us also mention that Calkin and Taylor \cite{CT} showed that the number of subsets of $[N]$ containing no solutions to $x + y + z = w$ is $O(2^{2N/3})$, an estimate which is basically sharp.\\[11pt]
It is natural to ask for estimates for $\mbox{SF}(\Gamma)$, the number of sum-free subsets of some finite abelian group $\Gamma$. When $\Gamma = \mathbb{Z}/p\mathbb{Z}$ ($p$ prime) this question is perhaps even more natural than the question of Cameron and Erd\H{o}s. It was first considered explicitly by Lev and Schoen \cite{LS}, who showed that 
\[ |\mbox{SF}(\mathbb{Z}/p\mathbb{Z})| \; \leq \; 2^{0.498p}.\]
Their result was improved by Ruzsa and the author \cite{GreenRuz1}, who obtained the estimate
\[ |\mbox{SF}(\mathbb{Z}/p\mathbb{Z})| \; \leq \; 2^{p/3 + o(p)}.\] This is tight except for the $o(p)$ term. For more general abelian groups $\Gamma$, work started by Lev, {\L}uczak and Schoen \cite{LLS} in the case $|\Gamma|$ even was continued by Ruzsa and the author \cite{GreenRuz2}, who obtained reasonably precise estimates for all abelian groups.\\[11pt]
The objective of the present paper is to prove the conjecture of Cameron and Erd\H{o}s. 
\begin{theoremenv}\label{mainthm}
The number of sum-free subsets of $[N]$ is asymptotically $c(N)2^{N/2}$, where $c(N)$ takes two different constant values according as $N$ is odd or even.
\end{theoremenv}
It is extremely likely that our methods extend to give, for example, much tighter bounds on $|\mbox{SF}(\mathbb{Z}/p\mathbb{Z})|$ but we do not pursue such matters here.\\[11pt]
\noindent\textbf{2. A strategy for counting sum-free sets.}
The purpose of this section is to outline the broad strategy that we will use to count sum-free subsets of $[N]$. Our method falls conveniently into two parts, which are dealt with in detail in the two sections immediately following this one. We have tried to make these sections as independent as possible. Our strategy, then, is as follows.\\[11pt]
Part I. We find some family $\mathcal{F}$ of subsets of $[N]$ with the following properties. Firstly, each $A \in \mathcal{F}$ is \textit{almost sum-free}, meaning that the number of \textit{additive triples} (triples with $x + y = z$) in $A$ is $o(N^2)$. Secondly, $\mathcal{F}$ does not contain too many sets; in fact, $|\mathcal{F}| = 2^{o(N)}$. Finally, every sum-free subset of $[N]$ is contained in some member of $\mathcal{F}$. \\[11pt]
Part II. Given $A \in \mbox{SF}(N)$, we consider some set $A' \in \mathcal{F}$ with $A \subseteq A'$. As $\mathcal{F}$ is so small, the number of $A$ for which $|A'| \leq \left(\frac{1}{2} - \frac{1}{120}\right)N$ is $o(2^{N/2})$. If, however, $|A'| \geq \left(\frac{1}{2} - \frac{1}{120}\right)N$ then it is possible to say something about the structure of $A'$, and hence about the structure of almost all $A \in \mbox{SF}(N)$.\\[11pt]
What we will actually show is that almost all $A \in \mbox{SF}(N)$ consist either entirely of odd numbers, or else are contained in the interval $\{\lceil (N + 1)/3\rceil,\dots,N\}$. The author was delighted to discover that, in their original paper \cite{CE1}, Cameron and Erd\H{o}s gave an elegant argument leading to an estimate for the number of sum-free subsets of $\{\lceil (N + 1)/3\rceil,\dots,N\}$. This argument, together with our work in the present paper, constitutes an affirmative solution to Conjecture \ref{conj1}.\\[11pt]
\noindent\textbf{3. Construction of $\mathcal{F}$. Granularizations.} 
In this section we complete Part I of the program outlined in \S 2 by constructing the family $\mathcal{F}$. This was basically achieved in the paper of Ruzsa and the author \cite{GreenRuz1}. Since it is not quite a trivial matter to isolate results from that paper in the form that we need them, we repeat some of the material from \cite{GreenRuz1} here.\\[11pt]
We begin with a small amount of notation concerning Fourier transforms. We will be working on the group $G = \mathbb{Z}/p\mathbb{Z}$, where $p$ is a prime. If $f : G \rightarrow \mathbb{C}$ is a function and if $r \in G$ then we define the Fourier transform $\hat{f}$ by
\[ \hat{f}(r) \; = \; \sum_{x \in G} f(x)e(rx/p)\] where, as usual, $e(\theta) = e^{2\pi i \theta}$. If $f,g$ are two functions then we define their convolution $f \ast g$ by \[(f \ast g)(x) \; = \; \sum_{y \in G} f(y)g(x - y).\] Observe that $(f \ast g)\hat{\;}(r) = \hat{f}(r)\hat{g}(r)$. Finally, we remark that if $A \subseteq G$ then we will identify $A$ with its characteristic function: that is, we write $A(x) = 1$ if $x \in A$ and $A(x) = 0$ otherwise. Observe that if $A,B \subseteq G$ are two sets then $(A \ast B)(x)$ is the number of representations of $x$ as $a + b$ with $a \in A$, $b \in B$.\\[11pt]
Let $p \in [2N,4N]$ be a prime, let $M$ be a positive integer, and let $d \in (\mathbb{Z}/p\mathbb{Z})^{*}$. Let $A \subseteq [N]$ be a set, and regard $A$ as a subset of $\mathbb{Z}/p\mathbb{Z}$ in the obvious manner. Suppose that $|A| = \alpha p$. Consider a partition of $\mathbb{Z}/p\mathbb{Z}$ into arithmetic progressions $I_i$, $i \in \mathbb{Z}/M\mathbb{Z}$, of common difference $d$ defined by
\begin{equation}\label{decomp} I_i \; = \; \left\{\lambda d \; : \; \frac{\overline{i}p}{M} \leq \lambda < \frac{\overline{(i+1)}p}{M}\right\},\end{equation} where $\overline{i}$ denotes the least positive residue of $i$. Each of these progressions has length either $L$ or $L-1$, where $L = \lceil p/M\rceil$. Let $\epsilon_1 > 0$ be a real number, and let \[ T \; = \; \{i \in \mathbb{Z}/M\mathbb{Z} \; | \; |A \cap I_i| \geq \epsilon_1|I_i|\}.\] Finally, define the \textit{granularization} $A'$ of $A$ (with respect to the length $d$ and the parameter $\epsilon_1$) by \[A' \; = \; \bigcup_{i \in T} I_i.\] It is easy to see that we have 
\begin{equation}\label{eq101} |A' \setminus A| \; \leq \; \epsilon_1 p.\end{equation}
One of the key results of \cite{GreenRuz1} is that, provided $d$ has a certain property (for which we will use the term ``good length'') the set $A'$ retains some of the additive features of $A$. In fact, we will be able to show that $A'$ is almost sum-free.\\[11pt]
Let us now say what we mean by the statement ``$d$ is a good length''. Let $\epsilon_2,\epsilon_3 > 0$ be two further real numbers and set $\delta  =\frac{1}{16}\epsilon_1^2\epsilon_2\epsilon_3^{1/2}\alpha ^{-1/2}$.
     Let $R$, $|R| = k$, be the set of all $r \neq  0$ for which $|\hat{A}(r) |\geq \delta  p
$. We say that $d$ is a \textit{good length} for $A$ (with respect to the parameters $\epsilon_1,\epsilon_2,\epsilon_3$) if
\begin{equation}\label{goodlength} \|dr/p\| \; \leq \;  {1\over 4L}\left( \frac{\delta p}{ |\hat{A}(r) |}\right )^{1/2}\end{equation}
for all $r\in R$. The following proposition was (essentially) the main result of \cite{GreenRuz1}. It clarifies the r\^ole of $\epsilon_2$ and $\epsilon_3$, which have not so far featured.
\begin{propositionenv}\label{prop111}
Suppose that
$d$ is a good length for $A$. Then the granularization $A'$ has the property that $A + A$ contains all $x$ for which $A' \ast A'(x) \geq \epsilon_2
p$, with at most $\epsilon_3 p$ exceptions.
\end{propositionenv}
\proof We claim that if $d$ is a good length then the function
\begin{equation}\label{eq766} g(x) = \frac{1}{2L - 1}\sum_{j = -(L-1)}^{L-1} e(jdx/p)\end{equation} satisfies
\begin{equation}\label{r1}
      |\hat{A}(x) | |1 - g(x)^2 |\;\leq \;\delta  p
    \end{equation}
     for all $x$. This automatically holds for $x=0$, as $g(x)=1$, and also
whenever $|\hat{A}(x)|\leq\delta p$, since $g(x)\in  [-1, 1]$. 
For any $x \in \mathbb{Z}/p\mathbb{Z}$ we may estimate $1-g(x)$ as follows. Writing $\Vert t\Vert$ for the distance of $t$ from
the nearest integer we have the inequality
    $1-\cos 2\pi t\leq  2\pi^2\| t\|^2$. It follows that
    \begin{eqnarray}
    \nonumber 1-g(x) & = & \frac{2}{2L-1}\sum_{j = 1}^{L-1}\left( 1 -\cos \frac{2\pi jdx}{ p}\right)\\
    \nonumber &\leq  & \frac{4\pi^2}{2L-1}\sum_{j = 1}^{L-1}\left\|  \frac{jdx}{p}\right\|^2\\
    \nonumber &\leq  & \frac{4\pi^2}{2L-1} \left\| \frac{dx}{p}\right\|^2\sum_{j = 1}^{ L-1}  j^2\\
    \label{r2} &\leq  & \frac{2\pi^2L^2}{3}\left\|  \frac{dx}{p}\right\|^2.
    \end{eqnarray}
     Hence
     \begin{eqnarray*}
      |\hat{A}(x) | |1 - g(x)^2 | & \leq  & 2 |\hat{A}(x) | |1 - g(x) | \\
      & \leq  & 14 L^2\| dx/p\|^2  |\hat{A}(x) |
    \end{eqnarray*}
 It is now easy to see that $d$ being a good length is exactly the property required to make \eqref{r1} hold.\\[11pt]
 Now to establish the proposition we define a function $a_1$ by
    \[a_1(n) \; = \;\frac{1}{|P|}(A\ast P)(n)\; =\;\frac{1}{ |P |} |A\cap(P
+ n) |,\]
     where $P =\{-(L-1)d,\dots,0,d,2d,\dots, (L-1)d\}$.
     Observe that $\hat{a}_1(x) = \hat{A}(x)g(x)$. Thus we have, by two applications of Parseval's identity, that
    \begin{eqnarray}
    \sum_n\left  |(A \ast A)(n) - (a_1 \ast a_1)(n)\right  |^2\nonumber & = & p^{-1}\sum_x\left  | \hat{A}(x)^2 - \hat{a}
_1(x)^2\right  |^2\\\nonumber & = & p^{-1}\sum_x  |\hat{A}(x) |^4\left(1 - g(x)^2\right )^2\\
&\nonumber\leq  & p^{-1}\left(\sup_x  |\hat{A}(x) | |1 - g(x)^2 |\right )^2\sum  _x  |\hat{A}(x) |^2\\
\label{inamin}& = &\alpha  p\left(\sup_x  |\hat{A}(x) | |1 - g(x)^2 |\right )^2.\end{eqnarray}
(\ref{r1}) therefore implies that
    \begin{equation}\label{eq377}
    \sum_n\left  |(A \ast A)(n) - (a_1 \ast a_1)(n)\right  |^2\;\leq \;\alpha \delta ^2 p^3.
    \end{equation}
     Now if $n\in A'$ then there is a progression of common difference $d$ and length $L$ containing $n$ which contains at least $\epsilon_1 L/2$
points of $A$. This progression is contained in $[n-(L-1)d,\dots, n+(L-1)d]$. Hence $a_1(n)$ is certainly at least $\epsilon_1/4$, and so $a
_1(n)\geq \epsilon_1 A(n)/4$ for all values of $n$. It follows immediately that $(a_1 \ast a_1)(n)\geq \epsilon_1^2 (A' \ast A')(n)/16$ for all $n$, and hence that if $A' \ast A'(n)\geq \epsilon_2 p$
then $a_1 \ast a_1(n)\geq
    \epsilon_1^2\epsilon_2 p/16$. We are to show that there are not many points $n$ for
which this is true whilst $A \ast A(n) = 0$. Letting $B$ denote the set of these
``bad'' points, observe that $n\in B$ implies that
     \[ \left  |(A \ast A)(n) - (a_1 \ast a_1)(n)\right  |^2 \; \geq  \;
\frac{\epsilon_1^4\epsilon_2^2p^2}{256}.\]
     Substituting into (\ref{eq377}) gives the bound
    \[  |B |\;\leq \;\frac{256\alpha \delta ^2}{\epsilon_1^4\epsilon_2^2} p\;\le \;\epsilon_3 p\]
    (this explains our choice of $\delta $). \endproof\\[11pt]
We defer for a while the issue of whether there are any good lengths. Our next result says that the conclusion of Proposition \ref{prop111} is enough to guarantee that if $A$ is sum-free then $A'$ is almost sum-free.
\begin{propositionenv}\label{almostsum-free}
Suppose that $A$ is sum-free. Let $\epsilon > 0$, set $\epsilon_1 = \epsilon$, $\epsilon_2 = \epsilon^2/144$, $\epsilon_3 = \epsilon^2/80$ and let $A'$ be the granularization of $A$ with respect to some good length $d$. Then $A'$ contains at most $\epsilon p^2$ triples $(x,y,z)$ with $x + y = z$.
\end{propositionenv}
\proof 
The choice of $p$ (that is, $p \geq 2N$) guarantees that $A$ is sum-free when considered as a subset of $\mathbb{Z}/p\mathbb{Z}$. Suppose without loss of generality that $d = 1$, and suppose for a contradiction that the proposition is false. Recall the notation we set up at the start of the section, particularly the definitions of the intervals $I_i$ and the set $T \subseteq \mathbb{Z}/m\mathbb{Z}$. We begin by claiming that there are at least $\epsilon M^2/4$ triples $(i,j,k) \in T^3$ for which $i + j = k$ or $k + 1$. Indeed note that if $x + y = z$ and if $x \in I_i$, $y \in I_j$ and $z \in I_k$ then $i + j = k$ or $k + 1$. However for a fixed triple $(i,j,k)$ with this property there are at most $4p^2/M^2$ triples $(x,y,z)$, so our claim follows from a simple double count.
For definiteness suppose that there are at least $\epsilon M^2/8$ triples $(i,j,k) \in T^3$ with $i + j = k$ (the argument when there are many triples with $i + j = k + 1$ is very similar). Let $K$ be the set of all $k \in T$ for which $T \ast T(k) \geq \epsilon M/16$ so that, by an easy averaging argument, we have $|K| \geq \epsilon M /16$.

Suppose that $i + j = k$ with $i,j,k \in T$, and suppose that $z$ lies in the middle $(1 - \epsilon/2)$ of $I_k$. Then the number of representations of $z$ as $x + y$ with $x \in I_i,y\in I_j$ is at least $\epsilon p/8M$. Therefore if $k \in K$ we have $A' \ast A'(z) \geq \epsilon^2 p/144$. Note, however, that since $k \in T$ the middle $(1 - \epsilon/2)$ of $I_k$ contains at least $\epsilon p/4M$ points of $A$. 

We have now shown that there are at least $\epsilon^2 p/64$ elements $x \in A$ for which $A' \ast A'(x) \geq \epsilon^2 p/144$. By the property of $A'$ described in Proposition \ref{prop111} we see that $A + A$ contains an element of $A$, contrary to our assumption that $A$ is sum-free.\endproof\\[11pt]
We now look at the issue of finding a good length. 
\begin{propositionenv} Let $A \subseteq \mathbb{Z}/p\mathbb{Z}$ have cardinality $\alpha p$.
A good length for $A$ with parameters $\epsilon_1,\epsilon_2,\epsilon_3$ exists if \begin{equation}\label{smiley} p\; >\;
    ( 4L)^{256\alpha^{2 }\epsilon_1^{-4}\epsilon_2^{-2}\epsilon_{3}^{-1} }.\end{equation}
\end{propositionenv}
\proof It follows by a standard application of the pigeonhole principle that 
a $d$ satisfying \eqref{goodlength} exists if
    \begin{equation}\label{r3}
     p \; > \; (4L)^k\left(\prod_{r\in R} { |\hat{A}(r) |\over\delta p}\right )^{1/2} .
    \end{equation}
We claim that this inequality is a consequence of the hypothesis on $p,L,\epsilon_1,\epsilon_2$ and $\epsilon_3$ in the statement of the proposition.
Indeed, observe that Parseval's identity implies that
    \begin{equation} \label{r4}
    \sum  _{r \in R}  |\hat{A}(r) |^2 \; \le  \; \alpha  p^2,
    \end{equation}
from which the arithmetic-geometric mean inequality gives
     \[ \prod _{r \in R}  |\hat{A}(r) | \; \le  \; \left(\frac{\alpha
p^2}{k}\right )^{k/2}.\]
     It follows that the right side of (\ref{r3}) is at most
    \begin{equation}\label{r5}(4L\alpha^{1/4}\delta^{-1/2} k^{-1/4} )^k,\end{equation}
which is an increasing function of $k$ in the range $k < \left(\frac{256L^4}{e}\right)\frac{\alpha}{\delta^2}$. However another consequence of (\ref{r4}) is the inequality $k < \alpha/\delta^2$, and hence (\ref{r5}) is itself bounded above by 
$(4L)^{\alpha /\delta^2}$. Recalling our choice of $\delta$ confirms the claim, and hence there is a $d$ for which (\ref{goodlength}) holds. \endproof\\[11pt]
To get the conclusion of Proposition \ref{almostsum-free} we required $\epsilon_1 = \epsilon$, $\epsilon_2 = \epsilon^2/144$ and $\epsilon_3 = \epsilon^2/80$. It is an easy but slightly tedious task to check that if we put $\epsilon = (\log N)^{-1/11}$ and $M = \left\lfloor N\exp(-(\log N)^{1/12})\right\rfloor$ then, at least for $N$ sufficiently large, $A$ has at least one good length. For the remainder of the section we assume that the parameters $\epsilon$ and $M$ take these values.\\[11pt]
We are now in a position to define our family of sets $\mathcal{F}$. Take $\mathcal{F}$ to consist of all sets which can be formed in the following manner. For all $d \in (\mathbb{Z}/p\mathbb{Z})^{\ast}$ consider the decomposition \eqref{decomp} of $\mathbb{Z}/p\mathbb{Z}$ into progressions $I_i$ ($i \in \mathbb{Z}/m\mathbb{Z}$) with common difference $d$. Let $\mathcal{G}$ be the collection of sets which are unions of progressions $I_i$, for some $d$. Now throw away from $\mathcal{G}$ all those sets which have more than $\epsilon p^2$ additive triples, giving a new collection $\mathcal{H}$. Finally, let $\mathcal{F}$ consist of all subsets of $[N]$ which can be obtained by adding at most $\epsilon p$ elements to some $H \cap [N]$, $H \in \mathcal{H}$.\\[11pt]
This may seem complicated. It turns out, however, that we can rather easily establish the following rather clean proposition concerning $\mathcal{F}$ which contains all the information we need for subsequent sections.
\begin{propositionenv}\label{prop462}
The family $\mathcal{F}$ has the following properties:\\
\emph{(i)} Every member of $\mathcal{F}$ has at most $o(N^2)$ additive triples;\\
\emph{(ii)} If $A$ is sum-free then $A$ is contained in some member of $\mathcal{F}$;\\
\emph{(iii)} $|\mathcal{F}| \leq 2^{o(N)}$.
\end{propositionenv}
\proof (i) By definition every set in $\mathcal{H}$ has at most $\epsilon p^2$ additive triples, and thus the same is true of sets of the form $H \cap [N]$, $H \in \mathcal{H}$. By adding $\epsilon p$ elements to such an $H$, we cannot create more than $3\epsilon p^2$ new additive triples. The result follows from the fact that $p \leq 4N$.\\
(ii) Set $\epsilon_1 = \epsilon$, $\epsilon_2 = \epsilon^2/144$ and $\epsilon_3 = \epsilon^2/80$. Choose a good length $d$ for $A$ with respect to $\epsilon_1,\epsilon_2,\epsilon_3$, and consider the granularization $A'$ with respect to $d$ and $\epsilon_1$. By Proposition \ref{almostsum-free} this lies in $\mathcal{H}$, and the result follows from \eqref{eq101}. \\
(iii) There are $p-1$ choices for $d$, and then $2^M$ ways to pick elements of $\mathcal{G}$. Thus $|\mathcal{H}| \leq p2^M$, and so $|\mathcal{F}|$ is at most $p2^M$ times the number of subsets of $[N]$ of size at most $\epsilon N$. This is clearly $2^{o(N)}$.
\endproof\\[11pt]
\noindent\textbf{4. The structure of almost sum-free sets.} In this section we study large almost sum-free sets. The results may be regarded as ``almost'' versions of the results of Freiman \cite{Frei4}. Freiman's methods do not seem to generalise easily to almost sum-free sets, so we have been forced to devise our own arguments. We will need one further piece of notation. If $K$ is a positive real number and if $A \subseteq G$ is a subset of an abelian group, we will write $D(A,K)$ for the set of all $x \in G$ which have at least $K$ representations as $a - a'$ with $a,a' \in A$. We call this the set of $K$-popular differences of $A$.\\[11pt]
In this section the objects $I_i$, $M$ and $\epsilon$ are not the same as in the previous section.
\begin{propositionenv}
\label{prop15} Let $\epsilon = o(N)$ and
suppose that $A \subseteq [N]$ has at most $\epsilon N^2$ additive triples, and that $|A| = (\frac{1}{2} - \eta)N$ where $\eta  \leq 1/50$ \emph{(}$\eta$ is allowed to be negative\emph{)}. Then one of the following alternatives occurs:\\[11pt]
\emph{(i)} With the possible exception of at most $32\epsilon^{1/8}N$ elements, $A$ is contained in some interval of length $(\frac{1}{2} + 3\eta + 60\epsilon^{1/8})N$;\\
\emph{(ii)} At most $54\epsilon^{1/8} N$ elements of $A$ are even.
\end{propositionenv}
Throughout what follows we shall assume that $|A| = (\frac{1}{2} - \eta)N$, $\eta \leq 1/50$ and that $A$ has at most $\epsilon N^2$ additive triples.
\begin{lemmaenv}
We have \begin{equation}\label{star}\textstyle\frac{1}{2}\displaystyle|D(A,\epsilon^{1/2} N)| + |A| \; \leq \; N(1 + 2\epsilon^{1/2}).\end{equation}
\end{lemmaenv}
\proof We have $\left|\left(D(A,\epsilon^{1/2} N) \cap \mathbb{Z}_{>0}\right) \cap A\right| \leq \epsilon^{1/2} N$, or else $A$ would contain more than $\epsilon N^2$ additive triples. The result follows quickly from this and the observation that $d$ is $K$-popular if and only if $-d$ is.\endproof
\begin{lemmaenv}\label{lem2}
For all but at most $8\epsilon^{1/4}N$ values of $a$, at least $|A|-16\epsilon^{1/4}N$ of the differences $a - a'$ with $a' \in A$ lie in $D(A,32\epsilon^{1/2} N)$.
\end{lemmaenv}
\proof Consider the graph on vertex set $A$ in which $a$ is joined to $a'$ if $a - a'$ is not in $D(A,32\epsilon^{1/2} N)$. It has at most $64\epsilon^{1/2} N^2$ edges. The number of vertices with degree more than $16\epsilon^{1/4} N$ is thus at most $8\epsilon^{1/4} N$.\endproof\\[11pt]
The next lemma, which is an application of basic graph theory, is \cite{LLS}, \S 4, Proposition 1. We specialise the result to the case we need.
\begin{lemmaenv}[Lev,\L uczak, Schoen]\label{LLSlem} Let $S$ be a subset of an abelian group $\Gamma$, $|\Gamma| \leq N$. Suppose that $|D(S,8\epsilon^{1/2}N)| \leq 2|S| - 16\epsilon^{1/4}N$. Then there is a set $X \subseteq S$, $|S \setminus X| \leq 4\epsilon^{1/4}N$, with $X - X \subseteq D(S,8\epsilon^{1/2}N)$.
\end{lemmaenv}
Now partition $[N]$ into intervals $I_i$ such that the smallest element of $I_i$ is $\lfloor 2i\epsilon^{1/8} N\rfloor$. Let $j$ be minimal so that $I_j$ contains at least $9\epsilon^{1/4}N$ points of $A$. Then, by Lemma \ref{lem2}, there is some $m \in I_i$ such that at least $|A| - 16\epsilon^{1/4}N$ of the differences $a - m$, $a \in A$, lie in $D(A,32\epsilon^{1/2} N)$. Let $k$ be maximal so that $I_k$ contains at least $9\epsilon^{1/4}N$ points of $A$. Again, there is $M \in I_k$ so that at least $|A|-16\epsilon^{1/4}N$ of the differences $a - m$, $a \in A$, lie in $D(A,32\epsilon^{1/2} N)$. Clearly for at least $|A|-32\epsilon^{1/4}N$ values of $a$ \textit{both} $a - m$ and $a - M$ are popular. Furthermore
\[ |A \cap [1,m]| \; \leq \; \frac{9\epsilon^{1/4}N}{\epsilon^{1/8}} \; \leq \; 9\epsilon^{1/8}N,\] and a similar inequality holds for $|A \cap [M,N]|$. Thus there is a set $B \subseteq A$, with the following properties:
\begin{enumerate}
\item $B$ is contained in $\{m + 1,\dots,M\}$;
\item $|B| \geq |A| - 50\epsilon^{1/8}N$;
\item For all $b \in B$ the differences $b - m$ and $b - M$ are both in $D(A,32\epsilon^{1/2} N)$.
\end{enumerate}
Observe that the first and second of these points imply that $M - m > N/4$.
Now let $t = M - m$. Following Lev and Smeliansky \cite{LevSmel}, consider the projection map $\pi : \mathbb{Z} \rightarrow \mathbb{Z}/t\mathbb{Z}$. We note some simple facts about this map in a lemma.
\begin{lemmaenv}\label{facts}
\emph{(i)} $|\pi(A)|  \geq  |A|-50\epsilon^{1/8}N$;\\
\emph{(ii)} Let $\delta > 0$. If $d \in D(A,4\delta N)$ then $\pi(d) \in D(\pi(A),\delta N)$;\\
\emph{(iii)} If $d \in D(\pi(A),8\delta N)$ then some element of $\pi^{-1}(d)$ lies in $D(A,\delta N)$.
\end{lemmaenv}
\proof (i) Clearly $|\pi(A)| \geq |\pi(B)| = |B|$.\\
(ii) If $d = a - a'$ then $\pi(d) = \pi(a) - \pi(a')$. For different representations of $d$ as $a - a'$, certain of these representations of $\pi(d)$ may be the same. However, since $t > N/4$, no element of $\mathbb{Z}/t\mathbb{Z}$ has more than 4 preimages under $\pi$ which lie in $A$. The result follows.\\
(iii) If $\pi(d) = \pi(a_i) - \pi(a'_i)$ then $a_i - a_i' = d + \lambda_i t$ for some $\lambda_i \in \mathbb{Z}$. As $t > N/4$ and $-N < a_i - a_i' < N$ there are at most 8 possible values for $\lambda_i$. Thus for at least one value of $i$ there are $\delta N$ solutions to $a_i - a_i' = d + \lambda_i t$.\endproof\\[11pt]
It is immediate from part (iii) of this lemma that $|D(A,\delta N)| \geq |D(\pi(A),8\delta N)|$. However we can do better than this, since for several $d$ at least \textit{two} of the elements $\pi^{-1}(d)$ are popular differences. Indeed, for any $b \in B$ we have $b - m,b - M \in D(A,32\epsilon^{1/2} N)$, but $b - m \equiv b - M\pmod{t}$. Certainly $\pi(b - m) \in D(\pi(A),8\epsilon^{1/2}N)$ by Lemma \ref{facts}(ii). Thus, by Lemma \ref{facts}(iii), we have 
\begin{eqnarray}|D(A,\epsilon^{1/2} N)| & \geq & \nonumber |D(\pi(A),8\epsilon^{1/2}N)| + |B|\\ & \geq & |D(\pi(A),8\epsilon^{1/2}N)| + |A| - 50\epsilon^{1/8}N.\label{eq2} \end{eqnarray}
Combining this with \eqref{star} gives
\begin{equation}\label{eq3} |D(\pi(A),8\epsilon^{1/2} N)| \;\leq \;2N(1 + 30\epsilon^{1/8}) - 3|A|.\end{equation}
Now we must have $|D(\pi(A),8\epsilon^{1/2} N)| \leq 2|\pi(A)| - 16\epsilon^{1/4}N$, since otherwise \eqref{eq3} and Lemma \ref{facts}(i) would give $|A| \leq (\frac{2}{5} + 100 \epsilon^{1/8})N$, which is contrary to our assumption about $|A|$. Thus Lemma \ref{LLSlem} applies, and we may pass to a subset $X \subseteq \pi(A)$ with \begin{equation}\label{eq5}|X| \geq |A| - 54\epsilon^{1/8}N\end{equation} and \begin{equation}\label{eq44}X - X \subseteq D(\pi(A),8\epsilon^{1/2} N).\end{equation} We distinguish three further cases.\\[11pt]
Case 1. $|X| \geq t/2$. Then $X - X$ is all of $\mathbb{Z}/t\mathbb{Z}$, and so \eqref{eq3} and \eqref{eq44} yield
\begin{equation}\label{eq4} t \leq \left(\textstyle\frac{1}{2}\displaystyle + 3\eta + 60\epsilon^{1/8}\right)N.\end{equation} 
But we know that, with at most $18\epsilon^{1/8}N$ exceptions, the elements of $A$ lie in the interval $\{m+1,\dots,M\}$ which has length $t$. 
This is alternative (i) of Proposition \ref{prop15}. \\[11pt]
Case 2. $|X - X| \geq 2|X| - t/3$. Then \eqref{eq3},\eqref{eq5} and \eqref{eq44} give $|A| \leq (\frac{7}{15} + 40\epsilon^{1/8})N$, contrary to assumption.\\[11pt]
Case 3. $|X - X| < 2|X| - t/3$. Then, by Kneser's theorem on the addition of sets in abelian groups (see \cite{Nath}, Theorem 4.2), $X - X$ is a union of cosets of some subgroup $H \leq \mathbb{Z}/t\mathbb{Z}$ of index 2. Thus $t$ is even and $\pi^{-1}(X)$ consists of integers of just one parity. That is, either at least $|A| - 54\epsilon^{1/8}N$ elements of $A$ are odd, or else at least that many are even. The latter possibility is, however, easily excluded; any subset of $\{2,4,6,\dots,2\lfloor N/2\rfloor\}$ of cardinality at least $12N/25$ contains at least $N^2/100$ additive triples. This concludes the proof of Proposition \ref{prop15}.\\[11pt]
An immediate corollary of Propositions \ref{prop462} and \ref{prop15} is the following result of Alon \cite{Alon}, Calkin \cite{Calk} and Erd\H{o}s and Granville (unpublished).
\begin{propositionenv}[Alon,Calkin,Erd\H{o}s--Granville] \label{AC} $|\mbox{\emph{SF}}(N)| = 2^{N/2 + o(N)}$.
\end{propositionenv}
\proof It follows from Proposition \ref{prop15} that if $F \subseteq [N]$ has $o(N^2)$ triples then $|F| \leq (\frac{1}{2} + o(1))N$. The result now follows from Proposition \ref{prop462}.\endproof\\[11pt]
A much more important corollary for us will be the following description of almost all sum-free subsets of $[N]$.
\begin{corollaryenv}\label{cor4} With $o(2^{N/2})$ exceptions, all sum-free subsets of $[N]$ consist entirely of odd numbers, or else are contained in $\{\lceil (N+1)/3\rceil,\dots,N\}$.
\end{corollaryenv}
\proof Let $A \in \mbox{SF}(N)$, and let $F \in \mathcal{F}$ contain $A$. The number of $A$ for which $|F| \leq \left(\frac{1}{2} - \frac{1}{120}\right)N$ is certainly $o(2^{N/2})$, so suppose that $|F| \geq \left(\frac{1}{2} - \frac{1}{120}\right)N$. Proposition \ref{prop15} then applies. Suppose first of all that alternative (ii) of that proposition holds, so that $A$ contains $o(N)$ even numbers. Suppose that $A$ contains at least one even number, $t$ say. If $t < N/2$ then we may select $\lfloor N/8\rfloor$ disjoint pairs $(x,x+t)$ of odd numbers, and $A$ cannot contain both of the elements of any of them since it is sum-free. The number of choices for $A$ is thus no more than $2^{N/4 + o(N)}3^{N/8} = o(2^{N/2})$. If $t \geq N/2$ then a very similar argument applies with pairs $(x,t - x)$. Thus all but $o(2^{N/2})$ of the sum-free sets with $o(N)$ even numbers consist \textit{entirely} of odd numbers.

Now suppose that alternative (i) of Proposition \ref{prop15} holds. If 
\begin{equation}\label{886} \left| A \cap \left[\left(1 - \textstyle\frac{1}{120}\displaystyle\right)N,N\right]\right| \; \leq \; 32\epsilon^{1/8}N\end{equation} then, using the fact that $A \cap [1,\left(1 - \frac{1}{120}\right)N]$ is sum-free together with Theorem \ref{AC}, we see that there are just $o(2^{N/2})$ possibilities for $A$. Suppose, then, that \eqref{886} fails to hold. Since Proposition \ref{prop15}, (i), holds we infer that $A$ is contained in the interval $[\left(\frac{1}{2} - \frac{1}{30} - 256\epsilon^{1/8}\right)N,N]$ with the exception of at most $32\epsilon^{1/8}N$ elements. Suppose that $A$ contains some element $t \in \{1,\dots,\lfloor (N+1)/3\rfloor\}$. Then we may select $\lfloor N/12\rfloor - 4$ totally disjoint pairs $(x,x+t)$ with $x \geq N/2$, and $A$ can contain at most one element from each of them. This means that the number of choices for $A$ is no more than $3^{N/12}2^{\left(\frac{1}{3} + \frac{1}{30} + o(1)\right)N}$ which, it can be checked, is $o(2^{N/2})$.\endproof\\[11pt]
As we remarked in the introduction, Cameron and Erd\H{o}s \cite{CE1} addressed the issue of counting sum-free subsets of $\{\lceil (N+1)/3\rceil,\dots,N\}$. They discovered that the number of such sets is asymptotically $c(N)2^{N/2}$, where $c(N)$ takes two different constant values depending on whether $N$ is odd or even. Combining their result with the work of this paper, then, leads to Theorem \ref{mainthm}.\endproof\\[11pt]
\noindent\textbf{Concluding remarks.} The paper \cite{CE1} of Cameron and Erd\H{o}s can be hard to locate and so we have written up their argument and posted it on the web \cite{GreenWeb}. The author would like to thank Imre Ruzsa for the many conversations which led to the papers \cite{GreenRuz1,GreenRuz2} and which, naturally, have had a significant bearing on the present work.

\noindent Ben Green\\
Trinity College, Cambridge, England.\\
email: bjg23@hermes.cam.ac.uk

     \end{document}